\documentclass[12pt,reqno]{amsart}
\usepackage{fullpage}
\usepackage{amsmath,amssymb,amsfonts,amsthm,mathtools}
\usepackage{booktabs}
\usepackage[colorlinks=true,citecolor=blue,linkcolor=blue,urlcolor=blue]{hyperref}
\usepackage{times}
\usepackage[utf8]{inputenc}
\newtheorem{theorem}{Theorem}[section]
\newtheorem{proposition}[theorem]{Proposition}
\newtheorem{lemma}[theorem]{Lemma}

\theoremstyle{definition}

\newtheorem{example}[theorem]{Example}
\newtheorem{remark}[theorem]{Remark}
\newcommand{\F}{\mathbb F}
\newcommand{\rk}{\operatorname{rank}}
\newcommand{\row}{\operatorname{row}}
\newcommand{\one}{\mathbf 1}
\newcommand{\adj}{\operatorname{adj}}
\newcommand{\Res}{\operatorname{Res}}

\title{Two-distance sets over finite fields}
\author{J\'ozsef Solymosi}
\address{Department of Mathematics, University of British Columbia, Vancouver, Canada \& \'Obuda University, Budapest, Hungary}
\email{solymosi@math.ubc.ca}
\author{Chi Hoi Yip}
\address{Department of Mathematics, Hong Kong University of Science and Technology, Clear Water Bay, Hong Kong}
\email{machyip@ust.hk}

\begin{document}

\begin{abstract}
We study two-distance sets in standard $d$-dimensional quadratic spaces over finite fields. In
characteristic $3$, we construct sets attaining the full
Larman--Rogers--Seidel bound, showing that Blokhuis' Euclidean bound $\binom{d+2}{2}$ need not
hold over finite fields. We prove a rank-sensitive replacement which precisely
measures the failure of positive definiteness. We also show that the
Blokhuis bound $\binom{d+2}{2}$ is attained over a suitable finite
field exactly for each $d\neq6$; in the exceptional dimension $d=6$, the exact
maximum is $27$.
\end{abstract}
\maketitle

\section{Introduction}
A finite set \(X\subset\mathbb R^d\) is called an \(s\)-distance set if
exactly \(s\) Euclidean distances occur between distinct points of \(X\).  The
classical problem is to bound \(|X|\) in terms of \(d\).  For a two-distance set
\(X\subset\mathbb R^d\), a standard polynomial argument, implicit in the work
of Larman, Rogers and Seidel~\cite{LRS}, gives
\[
        |X|\le \binom{d+2}{2}+d+1,
\]
while Blokhuis~\cite{Blokhuis} sharpened this to
\begin{equation}\label{eq:blokhuis-real}
        |X|\le \binom{d+2}{2}.
\end{equation}
More generally, Bannai, Bannai and Stanton~\cite{BBS} proved that every
\(s\)-distance subset of \(\mathbb R^d\) has size at most
\(\binom{d+s}{s}\); Petrov and Pohoata~\cite{PetrovPohoata} later gave a
short proof of this bound.  For sets constrained to a sphere, sharper bounds
are available, and the problem is closely connected with spherical codes and
association schemes; see Delsarte--Goethals--Seidel~\cite{DGS},
Delsarte--Larman~\cite{DL}, and the later refinements in
\cite{Musin,BargYu,GlazyrinYu}.

These universal bounds are rarely attained, and even the two-distance problem
is sensitive to the dimension.  The edge midpoints of a regular \(d\)-simplex
form a two-distance set of size \(\binom{d+1}{2}\).  Lison\v{e}k~\cite{Lisonek}
constructed a \(45\)-point two-distance set in \(\mathbb R^8\), attaining the
Blokhuis bound \(\binom{10}{2}\); see also
\cite{BannaiBannaiLeeXiangYu} for a recent study of the coherent-configuration
structure of this example.  More recently, Ge, Koolen and
Munemasa~\cite{GKM} constructed a \(277\)-point two-distance set in
\(\mathbb R^{23}\), exceeding the simplex-midpoint size \(276\) but remaining
below the Blokhuis bound \(300\).  Thus, even over \(\mathbb R\), determining
the extremal size and describing the exceptional configurations are delicate
problems.

Few-distance problems have also been studied in other ambient spaces.  In the
Hamming and Johnson spaces they become, respectively, problems about codes
with few Hamming distances and constant-weight codes with few distances; see
\cite{BargMusin,BargGlazyrinKaoLaiTsengYu,MusinNozaki}.  Blokhuis and
Seidel~\cite{BlokhuisSeidel} studied few-distance sets in the
pseudo-Euclidean spaces \(\mathbb R^{p,q}\).  Lison\v{e}k~\cite{Lisonek} also
gave a Lorentzian counterpart to his Euclidean constructions: for every
\(d\ge10\), he constructed a set of unit vectors in \(\mathbb R^{d,1}\) with
two inner products whose size attains the corresponding Blokhuis bound.  These
variants illustrate that the extremal behavior depends strongly on the
ambient quadratic form, and in particular on its positivity properties.

Finite fields provide a natural setting in which to separate the algebraic
part of the Euclidean theory from the role of positivity.  Let \(q\) be odd
and equip \(\mathbb F_q^d\) with the standard quadratic form
\[
        Q(x)=x_1^2+\cdots+x_d^2,
        \qquad \operatorname{dist}^2(x,y)=Q(x-y).
\]
A set \(S\subseteq\mathbb F_q^d\) is a two-distance set if
\(Q(x-y)\), for distinct \(x,y\in S\), takes exactly two values.  In our
upper-bound results these values are assumed to be distinct and nonzero.
Related questions about unit distances over finite fields were studied by
Blokhuis, Brouwer and Wilbrink~\cite{BBW}.  The standard polynomial argument
remains valid over every field of characteristic different from \(2\), but
the positive-definiteness step in Blokhuis' proof fails over finite fields.
This leads naturally to two questions: what replaces the Blokhuis bound, and
how large can finite-field two-distance sets be?

Interestingly, the Euclidean Blokhuis bound can fail by
the largest amount permitted by the standard polynomial argument.
Proposition~\ref{prop:lrs-bound} records the following field-independent form
of the dimension count (implicit in Larman--Rogers--Seidel~\cite{LRS}): if
\(K\) has characteristic different from \(2\) and
\(S\subset K^d\) is a two-distance set with two distinct nonzero squared
distances, then
\[
        |S|\le \binom{d+2}{2}+d+1.
\]
Theorem~\ref{thm:sharp-lrs-char3} gives an explicit construction showing that equality holds for every
\(q=3^r\) and every \(d\ge2\) with \(d\equiv2\pmod3\). Thus the finite-field extremal problem can differ sharply from
the real one. Theorem~\ref{thm:refined-rank-sensitive} explains this phenomenon by replacing
positivity with a rank defect.  For such a set \(S\), let \(H\) be its affine
coordinate matrix.  Then
\[
        |S|\le \binom{d+2}{2}+\rk H-\rk(HH^T),
\]
and a sharper inequality is obtained by adjoining to \(H\) the row of squared
norms.  The defect \(\rk H-\rk(HH^T)\) equals
\(\dim(\row(H)\cap\ker H)\).  It vanishes over \(\mathbb R\), recovering
inequality~\eqref{eq:blokhuis-real}, whereas the characteristic \(3\) examples attain
both rank-sensitive bounds.

We then turn to constructions attaining the Blokhuis bound.
Proposition~\ref{prop:standard-equilateral-criterion} characterizes exactly
when the standard space contains a \((d+2)\)-point equilateral set with
nonzero common squared distance.  Writing \(p=\operatorname{char}\mathbb F_q\),
this occurs if and only if
\[
        p\mid d+2
        \qquad\text{and}\qquad
        \bigl(d\text{ is odd or }q\equiv1\pmod4\bigr).
\]
Under these hypotheses, Theorem~\ref{thm:k-subset-sums-finite-field}, applied
with \(k=2\), gives a two-distance set of size \(\binom{d+2}{2}\).  More
generally, for each fixed \(s\), choosing \(p>s\) and sufficiently large
\(d\) among the infinitely many dimensions satisfying the criterion, the case
\(k=s\) gives an \(s\)-distance set of size \(\binom{d+2}{s}\).  These sets
asymptotically attain the real Bannai--Bannai--Stanton bound as
\(d\to\infty\).  Proposition~\ref{prop:midpoint-refined-sharp} shows that the
refined rank-sensitive estimate is exact for the midpoint examples.

A second Blokhuis-sized family, inspired by Lison\v{e}k's extension of the
Johnson representation \cite{Lisonek}, is constructed in
Proposition~\ref{prop:johnson-coset}. It applies whenever
\(p=\operatorname{char}\mathbb F_q>3\) and \(d\equiv8\pmod p\). Combining
this family with the characteristic \(3\) and equilateral-midpoint
constructions, Theorem~\ref{thm:blokhuis} completely determines the
problem of attaining the Blokhuis bound. For every \(d\ne6\), some standard space
\(\mathbb F_q^d\) contains a two-distance set of size
\(\binom{d+2}{2}\); in dimension $6$, every such set has at most \(27\)
points, and equality occurs over \(\F_{11}\).
\medskip

The paper is organized as follows. Section~\ref{sec:rank-sensitive} proves the standard polynomial
bound and its rank-sensitive refinement. Section~\ref{sec:char3} gives the characteristic
\(3\) construction and shows that both bounds are sharp. Section~\ref{sec:subset} develops
the equilateral and subset-sum constructions, and Section~\ref{sec:johnson} gives the
Johnson--coset construction. Finally, Section~\ref{sec:blokhuis} combines the
three construction families and resolves the exceptional dimension-$6$ case.

\section{Polynomial upper bounds and rank-sensitive refinements}\label{sec:rank-sensitive}

Throughout this section, let \(K\) be a field with \(\operatorname{char}K\ne2\), and put
\(Q(x)=x_1^2+\cdots+x_d^2\).

The following standard dimension count is implicit in the polynomial argument of
Larman--Rogers--Seidel~\cite{LRS} and remains valid over every field of
characteristic different from \(2\).

\begin{proposition}
\label{prop:lrs-bound}
Let \(S\subset K^d\) be a two-distance set with two distinct nonzero squared distances
\(\alpha,\beta\in K^*\).  Then
\[
        |S|\le \binom{d+2}{2}+d+1.
\]
\end{proposition}

\begin{proof}
Write \(S=\{s_1,\ldots,s_m\}\).  For each \(s\in S\), define
\[
        F_s(x)=\bigl(Q(x-s)-\alpha\bigr)\bigl(Q(x-s)-\beta\bigr).
\]
The evaluation matrix \(\bigl(F_{s_i}(s_j)\bigr)_{i,j}\) is diagonal with
nonzero diagonal entries \(\alpha\beta\).  Hence the polynomials \(F_s\),
\(s\in S\), are linearly independent.

Since
\[
        Q(x-s)=Q(x)-2\langle x,s\rangle+Q(s),
\]
each \(F_s\) belongs to
\[
        P_{\le2}+Q(x)\operatorname{span}\{x_1,\ldots,x_d\}+KQ(x)^2,
\]
where \(P_{\le2}\) is the space of polynomials of total degree at most \(2\).
This space has dimension \(\binom{d+2}{2}+d+1\), so the claimed bound follows.
\end{proof}

The following result makes the preceding obstruction explicit.  It bounds the possible
intersection between a row space and the kernel of the coordinate matrix; over
\(\mathbb R\) this intersection vanishes by positive definiteness, while over finite fields
it is measured by a rank defect. 

\begin{theorem}
\label{thm:refined-rank-sensitive} 
Let \(S=\{s_1,\ldots,s_m\}\subset K^d\) be a two-distance set with two distinct nonzero
squared distances \(\alpha,\beta\).  Let
\[
        H=
        \begin{pmatrix}
        1&1&\cdots&1\\
        s_{1,1}&s_{2,1}&\cdots&s_{m,1}\\
        \vdots&\vdots&&\vdots\\
        s_{1,d}&s_{2,d}&\cdots&s_{m,d}
        \end{pmatrix},
\]
and let \(\widehat H\) be obtained from \(H\) by adjoining the row
\((Q(s_1),Q(s_2),\ldots,Q(s_m))\).  Then
\begin{equation}\label{eq:refined-rank-bound}
        |S|\le
        \binom{d+2}{2}-1+\rk \widehat H-\rk(H\widehat H^T).
\end{equation}
In particular,
\begin{equation}\label{eq:affine-rank-bound}
        |S|\le
        \binom{d+2}{2}+\rk H-\rk(HH^T).
\end{equation}
Consequently, \(|S|\le \binom{d+2}{2}+d+1-\rk(HH^T)\).
\end{theorem}

\begin{proof}
For \(s\in S\), put
\[
        F_s(x)=\bigl(Q(x-s)-\alpha\bigr)\bigl(Q(x-s)-\beta\bigr).
\]
As in the proof of Proposition~\ref{prop:lrs-bound}, these polynomials are
linearly independent.  Let
\[
        W=\operatorname{span}\{F_s:s\in S\},
        \qquad
        B=\operatorname{span}\{1,x_1,\ldots,x_d,Q(x)\}.
\]
Then \(\dim W=|S|\), \(\dim B=d+2\), and \(W+B\) is contained in
\[
        P_{\le2}+Q(x)\operatorname{span}\{x_1,\ldots,x_d\}+KQ(x)^2,
\]
which has dimension \(\binom{d+2}{2}+d+1\).  Therefore
\begin{equation}\label{eq:need-bound-intersection}
        |S|\le \binom{d+2}{2}-1+\dim(W\cap B).
\end{equation}

It remains to control \(W\cap B\).  Suppose that
\(g=\sum_{i=1}^m c_iF_{s_i}\in B\).  Since \(g\) has degree at most \(2\),
comparison of the homogeneous parts of degrees \(4\) and \(3\) gives
\[
        \sum_{i=1}^m c_i=0,
        \qquad
        \sum_{i=1}^m c_is_i=0.
\]
Equivalently,
\begin{equation}\label{eq:Hc-zero}
        Hc=0,
        \qquad c=(c_1,\ldots,c_m)^T.
\end{equation}
On the other hand, \(g(s_i)=\alpha\beta c_i\) for every \(i\).  Since
\(g\in B\), the vector \((g(s_1),\ldots,g(s_m))\) lies in
\(\row(\widehat H)\); as \(\alpha\beta\ne0\), so does \(c\).  Thus
\[
        c\in\row(\widehat H)\cap\ker H.
\]
The representation of \(g\) in the basis \(\{F_{s_i}\}\) is unique, so the
map \(g\mapsto c\) is injective.  Hence
\[
        \dim(W\cap B)
        \le \dim\bigl(\row(\widehat H)\cap\ker H\bigr)
        =\rk\widehat H-\rk(H\widehat H^T),
\]
where the last equality follows because the image of the restriction of \(H\)
to \(\row(\widehat H)\) is \(\operatorname{im}(H\widehat H^T)\).  Together
with inequality~\eqref{eq:need-bound-intersection}, this proves
inequality~\eqref{eq:refined-rank-bound}.

To prove inequality~\eqref{eq:affine-rank-bound}, put \(U=\row(H)\) and
\(V=\row(\widehat H)\).  Then \(U\subseteq V\) and \(\dim(V/U)\le1\), while
\[
        \rk H-\rk(HH^T)=\dim(U\cap\ker H),
        \qquad
        \rk\widehat H-\rk(H\widehat H^T)=\dim(V\cap\ker H).
\]
Consequently,
\[
        \dim(V\cap\ker H)\le \dim(U\cap\ker H)+1.
\]
Substitution into inequality~\eqref{eq:refined-rank-bound} gives
inequality~\eqref{eq:affine-rank-bound}.  Finally, \(\rk H\le d+1\), which yields the
last assertion.
\end{proof}

The two most important special cases are the real setting and the spherical setting.  In
the real case the defect term disappears by positive definiteness.  In the spherical case
one can either use the refined rank-sensitive theorem over \(\mathbb R\), or give a direct
polynomial proof which works over every field of characteristic different from \(2\).

\begin{remark}\label{rem:spherical-bound}
Over \(\mathbb R\), the defect in inequality~\eqref{eq:affine-rank-bound} vanishes.  Indeed,
for every real matrix \(H\), one has \(\rk(HH^T)=\rk H\), since
\(u^THH^Tu=\|H^Tu\|^2\) for all real vectors \(u\).  This identity is the
linear-algebraic form of the positive-definiteness step in Blokhuis' real proof.  Thus
Theorem~\ref{thm:refined-rank-sensitive} recovers Blokhuis' bound in
inequality~\eqref{eq:blokhuis-real}.

More generally, let \(K\) be any field with \(\operatorname{char}K\ne2\), and suppose that
\(S\) is spherical.  After translating, write \(Q(s)=r\) for every \(s\in S\).  Then the
extra row \((Q(s_1),\ldots,Q(s_m))\) is \(r\) times the constant row of \(H\), so
\(\row(\widehat H)=\row(H)\).  Over \(\mathbb R\), the preceding paragraph makes the
refined defect zero, and Theorem~\ref{thm:refined-rank-sensitive} gives the usual
Delsarte--Goethals--Seidel bound~\cite{DGS}
\[
        |S|\le \binom{d+2}{2}-1.
\]
The same bound holds over every such field, although over a finite field the row-space
equality alone need not make the rank defect vanish.  Indeed, for
\(F_s(x)=(Q(x-s)-\alpha)(Q(x-s)-\beta)\), the identity
\(Q(x-s)=2r-2\langle x,s\rangle\) on \(S\) shows that each restriction
\(F_s|_S\) is represented by a polynomial in \(P_{\le2}\).  These restrictions are linearly
independent, since their evaluation matrix is \(\alpha\beta I\).  Since the nonzero
polynomial \(Q-r\in P_{\le2}\) vanishes on \(S\), the restriction space of \(P_{\le2}\) to
\(S\) has dimension at most \(\binom{d+2}{2}-1\), proving the claim.
\end{remark}

\section{Sharpness in characteristic 3}\label{sec:char3}

Proposition~\ref{prop:lrs-bound} gives \(|S|\le \binom{d+2}{2}+d+1\) for two-distance sets in \(\F_q^d\) with two distinct nonzero squared distances.  In
characteristic \(3\), this bound is sharp in infinitely many dimensions.

\begin{theorem}
\label{thm:sharp-lrs-char3}
Let \(q=3^r\), and let \(d\ge2\) satisfy \(d\equiv2\pmod3\).  On \(\F_q^d\), let
\(Q(x)=x_1^2+\cdots+x_d^2\) and \(\mathbf 1=(1,\ldots,1)\).  For \(1\le i<j\le d\), define \(u_{ij}=\mathbf 1-2e_i-2e_j\), so that \(u_{ij}\) has entries \(-1\) in positions \(i,j\) and entries \(1\) elsewhere.
For \(1\le i\le d\), define \(v_i=-\mathbf 1+2e_i\), so that \(v_i\) has entry \(1\) in position \(i\) and entries \(-1\) elsewhere.  Finally put
\[
        S_d=
        \{0,\mathbf 1\}
        \cup
        \{\pm e_i:1\le i\le d\}
        \cup
        \{u_{ij}:1\le i<j\le d\}
        \cup
        \{v_i:1\le i\le d\}.
\]
Then \(S_d\) is a two-distance set with nonzero squared distances \(1,2\in\F_q\), and
\(|S_d|=\binom{d+2}{2}+d+1\).  Consequently the finite-field Larman--Rogers--Seidel bound is sharp in infinitely many dimensions.
\end{theorem}

\begin{proof}
Because \(S_d\subset\F_3^d\), for any \(x,y\in S_d\) the value \(Q(x-y)\)
is the Hamming distance between \(x\) and \(y\), reduced modulo \(3\).
We show that every nonzero Hamming distance occurring in \(S_d\) is congruent
to \(1\) or \(2\pmod3\).

The distance from \(0\) to a point \(\pm e_i\) is \(1\), while its distance
to any of the sign vectors \(\mathbf1,u_{ij},v_i\) is
\(d\equiv2\pmod3\).  Two distinct points of
\(\{\pm e_i:1\le i\le d\}\) have Hamming distance \(1\) or \(2\).
Moreover, a point \(\pm e_i\) and a sign vector differ in every coordinate
other than possibly the \(i\)-th, so their Hamming distance is \(d-1\) or
\(d\), again congruent to \(1\) or \(2\pmod3\).

It remains to compare two sign vectors.  Their Hamming distance is the size of
the symmetric difference of their sets of negative coordinates, and the
possibilities are
\[
\begin{array}{c|c}
\text{pair type} & \text{Hamming distance} \\ \hline
\mathbf1,\ u_{ij} & 2 \\
\mathbf1,\ v_i & d-1 \\
u_{ij},u_{k\ell} & 2\text{ or }4 \\
u_{ij},v_k & d-1 \text{ if } k\in\{i,j\},\quad d-3 \text{ otherwise} \\
v_i,v_j & 2.
\end{array}
\]
Since \(d\equiv2\pmod3\), every listed value is congruent to \(1\) or
\(2\pmod3\).  Hence all nonzero squared distances in \(S_d\) lie in
\(\{1,2\}\).  Both occur, since
\[
        Q(e_i)=1,
        \qquad
        Q(\mathbf1)=d\equiv2\pmod3.
\]
Thus \(S_d\) is a two-distance set.

Finally, all points in the defining union are distinct.  The vectors
\(\pm e_i\) contain zero coordinates, whereas the remaining nonzero vectors
are sign vectors.  Among the latter, \(\mathbf1\), the \(u_{ij}\), and the
\(v_i\) have respectively \(0\), \(2\), and \(d-1\) negative coordinates;
these numbers are distinct for \(d\ge2\) with \(d\equiv2\pmod3\).  Within
each family, the negative coordinates determine the vector.  Therefore
\[
        |S_d|=2+2d+\binom d2+d
        =\binom{d+2}{2}+d+1,
\]
which is the bound in Proposition~\ref{prop:lrs-bound}.
\end{proof}

The same characteristic \(3\) examples also explain the size of the defect terms in the
rank-sensitive bounds.

\begin{remark}
The same examples also show that the rank-sensitive bounds are sharp.  Let \(H\) be the affine
coordinate matrix of \(S_d\).  The rank-sensitive bound gives
\[
        |S_d|\le
        \binom{d+2}{2}+\rk H-\rk(HH^T)
        \le
        \binom{d+2}{2}+d+1,
\]
since \(H\) has \(d+1\) rows.  But Theorem~\ref{thm:sharp-lrs-char3} gives equality in the last
bound.  Hence equality holds throughout; in particular, \(\rk H=d+1\) and \(\rk(HH^T)=0\).  Similarly, if \(\widehat H\) is obtained from \(H\) by adjoining the row
\((Q(x))_{x\in S_d}\), then the refined rank-sensitive bound gives
\[
        |S_d|\le
        \binom{d+2}{2}-1+\rk\widehat H-\rk(H\widehat H^T)
        \le
        \binom{d+2}{2}+d+1,
\]
because \(\widehat H\) has \(d+2\) rows.  Again equality holds throughout, so \(\rk\widehat H=d+2\) and \(\rk(H\widehat H^T)=0\).  Thus both the affine rank-sensitive bound and the refined rank-sensitive bound are attained.
\end{remark}

There is also a one-dimensional sharp example in characteristic \(5\).

\begin{example}
Let \(q=5^r\).  Take \(d=1\) and \(S=\F_5\subseteq\F_q\).  For distinct \(x,y\in S\), we have
\((x-y)^2\in\{1,4\}\).  Thus \(S\) is a two-distance set with two nonzero squared distances and
\(|S|=5=\binom{d+2}{2}+d+1\).  So the Larman--Rogers--Seidel upper bound is tight in this case.  Furthermore, if
\[
        H=
        \begin{pmatrix}
        1&1&1&1&1\\
        0&1&2&3&4
        \end{pmatrix},
\]
then \(HH^T=0\), because
\[
        \sum_{x\in\F_5}1=\sum_{x\in\F_5}x=\sum_{x\in\F_5}x^2=0.
\]
Hence the rank-sensitive bound gives equality.
\end{example}

\section{Equilateral and subset-sum constructions}\label{sec:subset}

\subsection{Equilateral sets over finite fields}

We next record the equilateral configurations that produce the midpoint examples.  The first
lemma is the standard Gram-matrix computation.

\begin{lemma}
\label{lem:equilateral-gram}
Let \(K\) be a field with \(\operatorname{char}K\ne2\).  Suppose
\(P_0,\ldots,P_{n-1}\in K^d\) satisfy \(Q(P_i-P_j)=\Delta\in K^*\) for all \(i\ne j\).
After translating so that \(P_0=0\), put \(v_i=P_i\) for \(1\le i\le n-1\).  Then \(Q(v_i)=\Delta\) and \(\langle v_i,v_j\rangle=\Delta/2\) for \(i\ne j\).  Equivalently, the Gram matrix of \(v_1,\ldots,v_{n-1}\) is
\[
        \frac{\Delta}{2}\bigl(I_{n-1}+J_{n-1}\bigr),
\]
where \(J_{n-1}\) is the all-ones matrix.
\end{lemma}

\begin{proof}
Since \(P_0=0\), we have \(Q(v_i)=\Delta\).  For \(i\ne j\),
\[
        \Delta=Q(v_i-v_j)
        =Q(v_i)+Q(v_j)-2\langle v_i,v_j\rangle,
\]
so \(\langle v_i,v_j\rangle=\Delta/2\).  The stated Gram matrix follows.
\end{proof}

This Gram-matrix computation immediately gives the standard upper bound for equilateral
sets over fields of odd characteristic.

\begin{proposition}
\label{prop:equilateral-size-bound}
Let \(K\) have odd characteristic \(p\), and let
\(P_0,\ldots,P_{n-1}\in K^d\) be equilateral with nonzero common squared distance.  Then \(n\le d+1\) unless \(p\mid n\).  In the exceptional case \(p\mid n\), one has \(n\le d+2\).
\end{proposition}

\begin{proof}
By Lemma~\ref{lem:equilateral-gram}, after translating one point to the origin
the Gram matrix is a nonzero scalar multiple of \(I_{n-1}+J_{n-1}\).  The
kernel of \(I_m+J_m\) is trivial unless \(p\mid m+1\), in which case it is
spanned by the all-ones vector.  Thus its rank is \(m\) unless
\(p\mid m+1\), and \(m-1\) otherwise.  Taking \(m=n-1\), the Gram matrix has
rank \(n-1\) if \(p\nmid n\), and rank \(n-2\) if \(p\mid n\).  Since its
rank is at most \(d\), the result follows.
\end{proof}

The previous proposition shows that the only possible way to have \(d+2\) equilateral points in
dimension \(d\) is the modular case \(p\mid d+2\).  In the standard quadratic space, there is one
additional square-class condition.  We use the following standard classification of
nondegenerate quadratic forms over finite fields; see, for example,
Lidl--Niederreiter~\cite[Section~6.2]{LidlNiederreiter}.

\begin{lemma}
\label{lem:finite-field-quadratic-forms}
Let \(q\) be odd.  Two nondegenerate quadratic forms over \(\F_q\) of the same dimension are
isometric if and only if their determinants have the same square-class in
\(\F_q^*/(\F_q^*)^2\).  Equivalently, if \(G\in \mathrm{M}_d(\F_q)\) is symmetric and nonsingular,
then there exists \(A\in\operatorname{GL}_d(\F_q)\) such that \(G=A^TA\) if and only if \(\det G\) is a square in \(\F_q^*\).
\end{lemma}

\begin{proposition}
\label{prop:standard-equilateral-criterion}
Let \(q\) be odd, let \(p=\operatorname{char}\F_q\), and let \(d\ge1\).  The standard quadratic
space \((\F_q^d,x_1^2+\cdots+x_d^2)\) contains a \((d+2)\)-point equilateral set with nonzero common squared distance if and only if
\begin{equation}\label{eq:equilateral-criterion}
        p\mid d+2
        \qquad\text{and}\qquad
        \bigl(d\text{ is odd }\text{or }q\equiv1\pmod4\bigr).
\end{equation}
\end{proposition}

\begin{proof}
Suppose first that \(P_0,\ldots,P_{d+1}\in\F_q^d\) are equilateral with
common squared distance \(\Delta\ne0\).  Translate so that \(P_0=0\), write
\(v_i=P_i\), and put \(\lambda=\Delta/2\).  By
Lemma~\ref{lem:equilateral-gram}, the Gram matrix of
\(v_1,\ldots,v_{d+1}\) is \(\lambda(I_{d+1}+J_{d+1})\).  Its rank is at most
\(d\), whereas
\[
        \det(I_{d+1}+J_{d+1})=d+2.
\]
Hence \(p\mid d+2\).

The Gram matrix of \(v_1,\ldots,v_d\) is
\(G=\lambda(I_d+J_d)\).  Since \(d+1=-1\) in \(\F_q\),
\[
        \det G=\lambda^d(d+1)=-\lambda^d.
\]
Writing \(G=V^TV\), we see that \(\det G\) is a square.  If \(d\) is even,
then \(\lambda^d\) is a square, so \(-1\) is a square in \(\F_q\),
equivalently \(q\equiv1\pmod4\).  This proves the necessity of
condition~\eqref{eq:equilateral-criterion}.

Conversely, assume condition~\eqref{eq:equilateral-criterion}.  Choose
\(\lambda\in\F_q^*\) so that \(-\lambda^d\) is a square: take
\(\lambda=1\) when \(d\) is even, and \(\lambda=-1\) when \(d\) is odd.
Then \(G=\lambda(I_d+J_d)\) is nonsingular with square determinant.  By
Lemma~\ref{lem:finite-field-quadratic-forms}, there is an
\(A\in\operatorname{GL}_d(\F_q)\) such that \(G=A^TA\).  Let
\(v_1,\ldots,v_d\) be the columns of \(A\).  Then
\[
        Q(v_i)=2\lambda,
        \qquad
        \langle v_i,v_j\rangle=\lambda\quad(i\ne j).
\]
Set \(v_{d+1}=-(v_1+\cdots+v_d)\).  Since \(d+2=0\) in \(\F_q\),
\[
        \langle v_i,v_{d+1}\rangle
        =-\bigl(2\lambda+(d-1)\lambda\bigr)=\lambda
\]
for every \(i\), and
\[
        Q(v_{d+1})=2d\lambda+d(d-1)\lambda=2\lambda.
\]
Thus \(0,v_1,\ldots,v_{d+1}\) is equilateral with nonzero common squared
distance \(2\lambda\).
\end{proof}

\subsection{Subset sums and midpoint constructions}

The same equilateral seed gives larger few-distance configurations.  We use subset sums
rather than averages; multiplying all points by a nonzero scalar only rescales squared
distances and does not change the number of distances.

\begin{proposition}
\label{prop:k-subset-sums}
Let \(K\) be a field with \(\operatorname{char}K\ne2\).  Let
\(A=\{v_1,\ldots,v_n\}\subset K^d\) be equilateral with common squared distance
\(\Delta\ne0\).  For each \(k\)-subset \(I\subset[n]\), define \(u_I=\sum_{i\in I}v_i\).  Let \(1\le k\le n-1\).  Then, for distinct \(k\)-subsets \(I,J\subset[n]\),
\[
        Q(u_I-u_J)=|I\setminus J|\Delta.
\]
Consequently, if \(\ell=\min(k,n-k)\) and
\(\Delta,2\Delta,\ldots,\ell\Delta\) are nonzero and pairwise distinct in \(K\), then \(\{u_I:|I|=k\}\) is an \(\ell\)-distance set of size \(\binom nk\).  In particular, this holds if
\(\operatorname{char}K=0\) or \(\ell<\operatorname{char}K\).
\end{proposition}

\begin{proof}
Let \(r=|I\setminus J|=|J\setminus I|\). For each \(i\), let \(c_i=1\) if \(i\in I\setminus J\), let \(c_i=-1\) if
\(i\in J\setminus I\), and let \(c_i=0\) otherwise. Then \(\sum_i c_i=0\), \(\sum_i c_i^2=2r\), and
\(u_I-u_J=\sum_i c_iv_i\).  Expanding the quadratic form and using
\(\sum_i c_i=0\) gives
\[
        \sum_{i,j}c_ic_jQ(v_i-v_j)
        =-2Q\left(\sum_i c_iv_i\right).
\]
Since \(Q(v_i-v_j)=\Delta\) for \(i\ne j\), the left-hand side is
\[
        \Delta\sum_{i\ne j}c_ic_j
        =-\Delta\sum_i c_i^2.
\]
Therefore
\[
        Q(u_I-u_J)
        =\frac{\Delta}{2}\sum_i c_i^2
        =r\Delta.
\]

For distinct \(I,J\), one has \(1\le r\le\ell\).  Under the stated
hypothesis, the values \(r\Delta\) are nonzero and pairwise distinct, so the
subset sums are distinct and determine an \(\ell\)-distance set.  Conversely,
for every \(1\le r\le\ell\), there are \(k\)-subsets \(I,J\) with
\(|I\setminus J|=r\), so all \(\ell\) distances occur.
\end{proof}

Applying this construction to maximal finite-field equilateral sets gives the following
finite-field few-distance examples.  The edge-midpoint construction is the case \(k=2\), after
rescaling the subset sums by the nonzero factor \(1/2\).

\begin{theorem}
\label{thm:k-subset-sums-finite-field}
Let \(q\) be an odd prime power, \(p=\operatorname{char}\F_q\), and \(d\ge1\).  Suppose that
\[
        p\mid d+2
        \qquad\text{and}\qquad
        \bigl(d\text{ is odd }\text{or }q\equiv1\pmod4\bigr).
\]
Let \(1\le k\le d+1\), and put \(\ell=\min(k,d+2-k)\).  If \(\ell<p\), then the standard quadratic space
\((\F_q^d,x_1^2+\cdots+x_d^2)\) contains an \(\ell\)-distance set of size \(\binom{d+2}{k}\).
In particular, if \(d\ge2\), then taking \(k=2\) gives a two-distance set of size
\(\binom{d+2}{2}\), namely the edge midpoints of a \((d+2)\)-point equilateral set.
\end{theorem}

\begin{proof}
By Proposition~\ref{prop:standard-equilateral-criterion}, the standard space
contains a \((d+2)\)-point equilateral set with nonzero common squared
distance \(\Delta\).  Apply Proposition~\ref{prop:k-subset-sums} with
\(n=d+2\).  Since \(\ell<p\), the elements
\(\Delta,2\Delta,\ldots,\ell\Delta\) are nonzero and pairwise distinct in
\(\F_q\), giving the required \(\ell\)-distance set.

If \(k=2\) and \(d\ge2\), then \(\ell=2<p\), because \(p\) is odd and
\(p\mid d+2\).  Multiplying the two-subset sums by \(1/2\) gives the edge
midpoints.
\end{proof}

For the midpoint examples, the refined rank-sensitive bound is exact.  The next computation
explains why the one-dimensional defect in the refined bound is exactly what is needed.

\begin{proposition}\label{prop:midpoint-refined-sharp}
Let \(q\) be odd, \(p=\operatorname{char}\F_q\), and \(d\ge2\) with
\(p\mid d+2\). Suppose that \(P_0,\ldots,P_{d+1}\in\F_q^d\) are equilateral
with nonzero common squared distance \(\Delta\), and put
\[
        X=\left\{\frac{P_i+P_j}{2}:0\le i<j\le d+1\right\}.
\]
Then \(X\) is a two-distance set of size \(\binom{d+2}{2}\), with squared
distances \(\Delta/4\) and \(\Delta/2\). If \(H_X\) is its affine coordinate matrix and
\(\widehat H_X\) is obtained by adjoining the row \((Q(x))_{x\in X}\), then
\[
        \rk\widehat H_X-\rk(H_X\widehat H_X^T)=1.
\]
Consequently, Theorem~\ref{thm:refined-rank-sensitive} gives the sharp bound
$|X|\le\binom{d+2}{2}.$
\end{proposition}

\begin{proof}
Proposition~\ref{prop:k-subset-sums}, applied with \(k=2\) and followed by
scaling by \(1/2\), gives the first assertion, including the cardinality of
\(X\).

Put \(N=d+2\) and \(\lambda=\Delta/2\), translate so that \(P_0=0\), and
let \(A=(P_1\,\cdots\,P_{N-1})\). By
Lemma~\ref{lem:equilateral-gram},
\[
        A^TA=\lambda(I_{N-1}+J_{N-1}).
\]
Since \(N=0\) in \(\F_q\), this matrix has rank \(N-2=d\). Thus
\(\rk A=d\), \(\ker A=\langle\mathbf1_{N-1}\rangle\), and
\[
        T:=\sum_{a=0}^{N-1}P_a=A\mathbf1_{N-1}=0.
\]
Indeed, \(A^TA\mathbf1_{N-1}=0\), so \(A\mathbf1_{N-1}\) is orthogonal to
the column space of \(A\), which is all of \(\F_q^d\).

Set
\[
        M:=\sum_{a=0}^{N-1}P_aP_a^T=AA^T.
\]
Then
\[
        M^2=A(A^TA)A^T=\lambda M,
\]
because \(A\mathbf1_{N-1}=0\). Moreover, \(M\) is invertible: if \(My=0\),
then \(A^Ty=c\mathbf1_{N-1}\) for some \(c\), and
\[
        0=T^Ty=\mathbf1_{N-1}^TA^Ty=(N-1)c=-c.
\]
Hence \(A^Ty=0\), and therefore \(y=0\). It follows that $       M=\lambda I_d.$

The first two moments of \(X\) are
\[
        \sum_{x\in X}1=\binom N2=0,
        \qquad
        \sum_{x\in X}x=\frac{N-1}{2}T=0,
\]
while
\[
\begin{aligned}
        \sum_{x\in X}xx^T
        =\frac14\sum_{i<j}(P_i+P_j)(P_i+P_j)^T
        =\frac14\bigl((N-2)M+TT^T\bigr)
          =-\frac{\Delta}{4}I_d.
\end{aligned}
\]
Thus
\[
        H_XH_X^T=
        \begin{pmatrix}
        0&0\\
        0&-\frac{\Delta}{4}I_d
        \end{pmatrix}.
\]

Now let \(c=(c_x)_{x\in X}\in\row(\widehat H_X)\cap\ker H_X\), say
\[
        c_x=a+b\cdot x+\gamma Q(x).
\]
Then
\[
        \sum_{x\in X}c_x=0,
        \qquad
        \sum_{x\in X}c_xx=0.
\]
Since
\[
        \sum_{x\in X}Q(x)
        =\operatorname{tr}\left(\sum_{x\in X}xx^T\right)
        =-\frac{d\Delta}{4}
        =\frac{\Delta}{2}\ne0,
\]
the first equation gives \(\gamma=0\). The second then gives
\[
        -\frac{\Delta}{4}b=0,
\]
so \(b=0\). Hence \(c\) is constant. Conversely, the constant vector belongs
to \(\row(\widehat H_X)\cap\ker H_X\), since
\(\sum_{x\in X}1=\sum_{x\in X}x=0\). Therefore
\[
\begin{aligned}
        \rk\widehat H_X-\rk(H_X\widehat H_X^T)
        =\dim\bigl(\row(\widehat H_X)\cap\ker H_X\bigr)
        =1.
\end{aligned}
\]
Theorem~\ref{thm:refined-rank-sensitive} now yields
\[
        |X|\le\binom{d+2}{2}-1+1=\binom{d+2}{2}.\qedhere
\]
\end{proof}

Here is a concrete example in characteristic \(7\) and dimension \(5\).

\begin{example}
In the standard space \(\F_7^5\), take \(\lambda=3\), so that \(\Delta=2\lambda=6\).  The
points \(P_0=0\) and
\[
\begin{aligned}
 v_1&=(0,0,0,2,3),&
 v_2&=(0,0,1,1,5),&
 v_3&=(0,1,5,0,1),\\
 v_4&=(2,5,5,0,1),&
 v_5&=(2,6,5,5,0),&
 v_6&=(3,2,5,6,4)
\end{aligned}
\]
form an equilateral set: \(Q(v_i)=6\) and \(\langle v_i,v_j\rangle=3\) for \(i\ne j\).
Therefore the \(21\) edge midpoints form a two-distance set in \(\F_7^5\), with squared
distances \(6/4=5\) and \(6/2=3\).  Thus the construction gives a two-distance set of size
\(21=\binom72=\binom{d+2}{2}\).
\end{example}

\section{A Lison\v{e}k-type Johnson--coset construction}\label{sec:johnson}
Lison\v{e}k constructed a maximum two-distance set of size \(45\) in
\(\mathbb R^8\) by enlarging the usual representation of the Johnson graph
\(J(9,2)\)~\cite{Lisonek}.  Bannai, Sato and Shigezumi later studied maximal
\(m\)-distance sets containing Johnson graph representations~\cite{BannaiSatoShigezumi}.
The notation \(J(n,2)\) refers to the standard Johnson layer, namely the set of
all two-subsets of an \(n\)-element set, or equivalently the binary vectors of
weight two; this terminology is rooted in Johnson's work on coding
bounds~\cite{Johnson}.  We now give a finite-field analogue of the same idea when $d \equiv 8 \pmod p$.
The construction starts from the weight-two Johnson layer and adds one further
\(n\)-point coset. 

\begin{proposition}\label{prop:johnson-coset}
Let \(q\) be a prime power of characteristic \(p>3\), and let \(d\ge1\) satisfy
\(d\equiv8\pmod p\). Write \(Q(x)=x_1^2+\cdots+x_d^2\) on \(\F_q^d\), and
put \(\mathbf1=(1,\ldots,1)\). Then
\[
\begin{aligned}
X_d={}&\left\{\frac12\mathbf1+e_i+e_j:1\le i<j\le d\right\}
\cup\left\{\frac14\mathbf1+e_i:1\le i\le d\right\}\\
&\cup\left\{\frac34\mathbf1-e_i:1\le i\le d\right\}\cup\{\mathbf1\}
\end{aligned}
\]
is a two-distance set with squared distances \(2\) and \(4\), and
\[
        |X_d|=\binom{d+2}{2}.
\]
\end{proposition}

\begin{proof}
Put \(n=d+1\), so \(n\equiv9\pmod p\), and work in \(\F_q^n\) with
\(Q_n(z)=\sum_{i=1}^n z_i^2\). Let \(\lambda=3^{-1}\) and
\[
        T_0=\{e_i+e_j:1\le i<j\le n\},\qquad
        T_1=\{\lambda\mathbf1_n-e_i:1\le i\le n\}.
\]
Every point of \(T=T_0\cup T_1\) has coordinate sum \(2\), since
\(n\lambda-1=2\). Hence all differences of points of \(T\) lie in the
zero-sum hyperplane.

Two distinct points of \(T_0\) have squared distance \(2\) or \(4\), while
two distinct points of \(T_1\) have squared distance \(2\). If
\(v=e_a+e_b\in T_0\) and \(w=\lambda\mathbf1_n-e_c\in T_1\), then
\[
        Q_n(v-w)=
        \begin{cases}
        5-6\lambda+n\lambda^2=4,&c\in\{a,b\},\\
        3-6\lambda+n\lambda^2=2,&c\notin\{a,b\}.
        \end{cases}
\]
Thus \(T\) is a two-distance set with squared distances \(2\) and \(4\).

Let
\[
        \mathcal H=\left\{z\in\F_q^n:\sum_{i=1}^n z_i=2\right\}.
\]
If \(z\) is a difference of two points of \(\mathcal H\) and
\(u=(z_1,\ldots,z_d)^T\), then \(z_n=-\sum_{i=1}^d z_i\), so
\[
        Q_n(z)=u^T(I_d+J_d)u.
\]
Set \(A=I_d+\frac14J_d\). Since \(J_d^2=dJ_d\) and \(d\equiv8\pmod p\),
\[
        A^TA
        =I_d+\left(\frac12+\frac d{16}\right)J_d
        =I_d+J_d.
\]
Moreover, \(A\) is invertible, since it acts as the identity on
\(\mathbf1^\perp\) and as multiplication by \(1+d/4=3\) on
\(\langle\mathbf1\rangle\). Hence
\[
        \Phi:\mathcal H\longrightarrow\F_q^d,\qquad
        \Phi(z)=A(z_1,\ldots,z_d)^T,
\]
is injective and preserves squared distances.

A direct calculation gives
\[
\begin{aligned}
\Phi(e_i+e_j)&=\frac12\mathbf1+e_i+e_j,
&
\Phi(e_i+e_n)&=\frac14\mathbf1+e_i,\\
\Phi(\lambda\mathbf1_n-e_i)&=\frac34\mathbf1-e_i,
&
\Phi(\lambda\mathbf1_n-e_n)&=\mathbf1,
\end{aligned}
\]
with the natural ranges for \(i,j\). Thus \(\Phi(T)=X_d\), and
\[
        |X_d|=|T_0|+|T_1|
        =\binom n2+n
        =\binom{d+2}{2}. \qedhere
\]
\end{proof}

\section{Configurations attaining the Blokhuis bound and the exceptional case \texorpdfstring{\(d=6\)}{d=6}}\label{sec:blokhuis}
The preceding constructions show that the Blokhuis bound is attained in every
dimension except $6$. We now combine them and determine the exact maximum in
dimension $6$.

\begin{theorem}\label{thm:blokhuis}
The following statements hold.
\begin{enumerate}
\item[\textup{(i)}] For every integer \(d\ge1\) with \(d\ne6\), there exists an
odd prime power \(q\) such that the standard quadratic space $        \bigl(\mathbb F_q^d,\ x_1^2+\cdots+x_d^2\bigr)$ 
contains a two-distance set of size \(\binom{d+2}{2}\).
\item[\textup{(ii)}] Let \(q\) be odd. Every two-distance set in the standard
quadratic space $       \bigl(\F_q^6,\ x_1^2+\cdots+x_6^2\bigr)$ has at most \(27\)
points. Equality occurs for example when \(q=11\).
\end{enumerate}
\end{theorem}

Part~\textup{(i)} follows by combining the preceding constructions.
The proof of part~\textup{(ii)} is completed by the exact computation in
Section~\ref{subsec:d6-computation}.

\subsection{Constructions outside dimension $6$}

\begin{proof}[Proof of Theorem~\ref{thm:blokhuis}\textup{(i)}]
First suppose that \(d\ge9\).  If \(d\equiv2\pmod3\), then
Theorem~\ref{thm:sharp-lrs-char3} gives a two-distance set
\(S_d\subseteq\F_3^d\) of size \(\binom{d+2}{2}+d+1\).  Choose a subset of
size \(\binom{d+2}{2}\) containing \(0,e_1,\mathbf1\).  Both distances still
occur, since
\[
        Q(e_1)=1,
        \qquad
        Q(\mathbf1)=d\equiv2\pmod3.
\]

Assume now that \(d\not\equiv2\pmod3\).  If \(d-8\) has an odd prime divisor
\(p\), then \(p\ne3\), and hence \(p>3\).  Since \(d\equiv8\pmod p\),
Proposition~\ref{prop:johnson-coset} gives the required set in \(\F_p^d\).

It remains to consider the case in which \(d-8\) has no odd prime divisor.
Then \(d-8=2^a\) for some integer \(a\ge0\).  The integer
\[
        d+2=2^a+10
\]
is not a power of \(2\): it equals \(11\) when \(a=0\), equals \(12\) when
\(a=1\), and is congruent to \(2\pmod4\) when \(a\ge2\).  It therefore has
an odd prime divisor \(p\).  Take
\(q=p^2\).  Then \(p\mid d+2\) and \(q\equiv1\pmod4\), so
Theorem~\ref{thm:k-subset-sums-finite-field} with \(k=2\) gives a two-distance
set of size \(\binom{d+2}{2}\) in \(\F_{p^2}^d\).

It remains to check \(d\le8\).  The cases \(d=3,4,7\) follow from
Theorem~\ref{thm:k-subset-sums-finite-field} with the following choices:
\[
\begin{array}{c|c|c}
d&q&\text{reason}\\ \hline
3&5&5\mid d+2,\ d\text{ odd},\\
4&9&3\mid d+2,\ q\equiv1\pmod4,\\
7&3&3\mid d+2,\ d\text{ odd}.
\end{array}
\]
For \(d=1\), Proposition~\ref{prop:johnson-coset} applies in characteristic
\(7\), since \(1\equiv8\pmod7\).  Finally, for \(d=2,5,8\), apply
Theorem~\ref{thm:sharp-lrs-char3} and again take a subset of size
\(\binom{d+2}{2}\) containing \(0,e_1,\mathbf1\).  In each case both squared
distances \(1\) and \(2\) remain.  This covers every \(d\ne6\).
\end{proof}

\subsection{A 27-point configuration in dimension $6$}\label{subsec:d6-construction}

\begin{proposition}\label{prop:d6-construction}
The standard quadratic space \(\F_{11}^6\) contains a \(27\)-point
two-distance set with squared distances \(2\) and \(4\).
\end{proposition}

\begin{proof}
Work first in \(\F_{11}^8\), with standard basis \(e_1,\ldots,e_8\), and put
\[
        A_i=e_i+e_7,\qquad B_i=e_i+e_8\quad(1\le i\le6),
\]
\[
        C_{ij}=\frac12\one_8-e_i-e_j\quad(1\le i<j\le6).
\]
Let \(\mathcal X\) be the set of these \(6+6+15=27\) vectors. They are
distinct. Every vector in \(\mathcal X\) has norm \(2\), and for distinct
vectors \(x,y\in\mathcal X\) one has \(\langle x,y\rangle\in\{0,1\}\). Indeed,
\begin{align*}
\langle A_i,A_j\rangle&=\langle B_i,B_j\rangle=1 &&(i\ne j),\\
\langle A_i,B_j\rangle&=\delta_{ij},\\
\langle A_k,C_{ij}\rangle=\langle B_k,C_{ij}\rangle
&=\begin{cases}0,&k\in\{i,j\},\\1,&k\notin\{i,j\},\end{cases}\\
\langle C_{ij},C_{k\ell}\rangle&=|\{i,j\}\cap\{k,\ell\}|
&&\bigl(\{i,j\}\ne\{k,\ell\}\bigr).
\end{align*}
Consequently,
\[
        Q(x-y)=Q(x)+Q(y)-2\langle x,y\rangle\in\{2,4\},
\]
and both values occur.

The set \(\mathcal X\) lies in the affine $6$-space
\[
        \mathcal A=\left\{z=(z_1,\ldots,z_8)^T\in\F_{11}^8:
        \sum_{i=1}^6z_i=1,\ z_7+z_8=1\right\}.
\]
Its direction space has basis
\[
        f_i=e_i-e_6\quad(1\le i\le5),\qquad g=e_7-e_8,
\]
with Gram matrix
\[
        G=\begin{pmatrix}I_5+J_5&0\\0&2\end{pmatrix}.
\]
Since
\[
        \det G=2\det(I_5+J_5)=2\cdot6=1\quad\text{in }\F_{11},
\]
Lemma~\ref{lem:finite-field-quadratic-forms} identifies the direction space
of \(\mathcal A\) with the standard quadratic space \(\F_{11}^6\). Translating
\(\mathcal X\) by the negative of any point of \(\mathcal A\), and then applying
this isometry, gives the required configuration.
\end{proof}

\subsection{Reduction to a finite graph problem}\label{subsec:d6-reduction}

Suppose for contradiction that
\[
        X=\{x_1,\ldots,x_{28}\}\subset K^6
\]
is a two-distance set over a finite field \(K\) of odd characteristic.
Throughout this subsection, vectors are column vectors. Unless another field
is specified, all linear algebra is over \(K\), and every \(0\)-\(1\) matrix is
viewed over \(K\). Let the squared
distances be distinct elements \(\alpha,\beta\in K^*\). After rescaling the
quadratic form by \((\alpha-\beta)^{-1}\), the two distances are
\[
        \lambda+1,\qquad \lambda,
        \qquad
        \lambda=\frac{\beta}{\alpha-\beta},
\]
where \(\lambda\ne0,-1\). We continue to denote the rescaled form by \(Q\)
and its associated bilinear form by \(\langle\ ,\ \rangle\).

Put
\[
        P=(x_1\ \cdots\ x_{28})\in K^{6\times28},
        \qquad
        \nu=(Q(x_1),\ldots,Q(x_{28}))^T\in K^{28}.
\]
The points are the columns of \(P\), and the columns of \(H\) and
\(\widehat H\) below are indexed by \(X\). The set \(X\) affinely spans
\(K^6\); otherwise, applying the proof of Proposition~\ref{prop:lrs-bound}
to the restriction of \(Q\) to the affine span of \(X\), whose dimension is
at most \(5\), gives
\[
        |X|\le\binom72+6=27.
\]
Moreover, \(X\) is not spherical by Remark~\ref{rem:spherical-bound}.
Therefore the affine coordinate matrix
\[
        H=\begin{pmatrix}\one_{28}^T\\ P\end{pmatrix}
\]
has rank \(7\), and the augmented coordinate matrix
\begin{equation}\label{eq:d6-Hhat}
        \widehat H=
        \begin{pmatrix}\one_{28}^T\\ P\\ \nu^T\end{pmatrix}
        \in K^{8\times28}
\end{equation}
has rank \(8\). Indeed, if \(\nu^T\in\row(H)\), then \(Q(x)=a+\ell(x)\) on \(X\) for some
linear functional \(\ell\), and completing the square would make \(X\)
spherical.

Let \(G\) be the graph on vertex set \(\{1,\ldots,28\}\) in which distinct
vertices \(i\) and \(j\) are adjacent if and only if
\(Q(x_i-x_j)=\lambda+1\). Let \(A\) be its adjacency matrix, and put
\[
        M=A-\lambda I_{28}.
\]
Then, for all \(1\le i,j\le28\),
\begin{equation}\label{eq:d6-Mentries}
        M_{ij}=Q(x_i-x_j)-\lambda.
\end{equation}

\begin{lemma}\label{lem:d6-kernel}
We have \(\rk M=8\). Moreover, there is a vector \(y\in K^{28}\) such that
\begin{equation}\label{eq:d6-affine-vector}
        My=\one_{28},
        \qquad
        \one_{28}^Ty=0.
\end{equation}
\end{lemma}

\begin{proof}
Let \(c=(c_1,\ldots,c_{28})^T\in K^{28}\), and set
\[
        s=\one_{28}^Tc,
        \qquad
        u=Pc,
        \qquad
        r=\nu^Tc.
\]
Thus
\[
        \widehat Hc=\begin{pmatrix}s\\u\\r\end{pmatrix}.
\]
Expanding equation~\eqref{eq:d6-Mentries} gives
\begin{equation}\label{eq:d6-moment-identity}
        (Mc)_i=(Q(x_i)-\lambda)s-2\langle x_i,u\rangle+r.
\end{equation}

We first prove
\[
        \ker M=\ker\widehat H.
\]
If \(\widehat Hc=0\), then \(s=u=r=0\), and
identity~\eqref{eq:d6-moment-identity} gives \(Mc=0\). Conversely, suppose
that \(Mc=0\). If \(s\ne0\), then
equation~\eqref{eq:d6-moment-identity} yields
\[
        Q\left(x_i-\frac{u}{s}\right)
        =\lambda-\frac{r}{s}+Q\left(\frac{u}{s}\right)
        \qquad(1\le i\le28),
\]
so \(X\) is spherical, a contradiction. Hence \(s=0\), and
equation~\eqref{eq:d6-moment-identity} becomes
\[
        2\langle x_i,u\rangle=r
        \qquad(1\le i\le28).
\]
The affine-linear function \(x\mapsto2\langle x,u\rangle-r\) therefore
vanishes on \(X\). Since \(X\) affinely spans \(K^6\), its linear part
vanishes, and the nondegeneracy of the bilinear form gives \(u=0\); then
\(r=0\). Thus \(\ker M=\ker\widehat H\), and hence \(\rk M=8\).

Since \(\widehat H\) has full row rank, choose \(y\in K^{28}\) such that
\[
        \widehat Hy=
        \begin{pmatrix}0\\ \mathbf0_6\\ 1\end{pmatrix}.
\]
Identity~\eqref{eq:d6-moment-identity} then gives
equation~\eqref{eq:d6-affine-vector}.
\end{proof}

\begin{proposition}\label{prop:d6-reduction}
Let \(p=\operatorname{char}K\). Then \(\lambda\in\F_p\). Moreover,
there exist a simple graph \(\Gamma\) on eight vertices, with adjacency
matrix \(C\in M_8(\F_p)\), and $20$ distinct vectors
\(b_1,\ldots,b_{20}\in\{0,1\}^8\subset\F_p^8\) such that
\[
        B=C-\lambda I_8
\]
is nonsingular and, with \(R=B^{-1}\),
\begin{align}
b_i^TRb_i&=-\lambda, \label{eq:d6-norm}\\
b_i^TRb_j&\in\{0,1\}\quad(i\ne j), \label{eq:d6-compat}\\
b_i^TR\one_8&=1, \label{eq:d6-main}\\
\one_8^TR\one_8&=0. \label{eq:d6-isotropic}
\end{align}
\end{proposition}

\begin{proof}
Write \(K=\F_{p^e}\). The adjacency matrix \(A\) has entries in
\(\F_p\). By Lemma~\ref{lem:d6-kernel}, the nullity of \(A-\lambda I_{28}\) is
\(20\), so \(\lambda\) is an eigenvalue of the \(\F_p\)-matrix \(A\)
with geometric multiplicity \(20\). If \(f=[\F_p(\lambda):\F_p]\), then
each Frobenius conjugate of \(\lambda\) is also an eigenvalue of \(A\) with
geometric, and hence algebraic, multiplicity at least \(20\). Since the
characteristic polynomial of \(A\) has degree \(28\), we have \(20f\le28\).
Thus \(f=1\), and \(\lambda\in\F_p\).

Thus \(M=A-\lambda I_{28}\) belongs to \(M_{28}(\F_p)\) and has rank \(8\)
over \(\F_p\). The block matrices below are therefore taken over \(\F_p\);
when necessary, we compare them with \(\widehat H\) and the vector \(y\) after
extending scalars to \(K\). Since \(M\) is symmetric of rank \(8\), it has a nonsingular principal submatrix of
order \(8\). After relabeling the vertices,
write
\[
        M=
        \begin{pmatrix}B&N\\N^T&E\end{pmatrix},
        \qquad
        B\in\F_p^{8\times8},\quad N\in\F_p^{8\times20},\quad
        E\in\F_p^{20\times20},
\]
and put \(R=B^{-1}\). The first eight vertices induce the graph \(\Gamma\), so
\(B=C-\lambda I_8\). Write
\[
        N=(b_1\ \cdots\ b_{20}),
        \qquad b_i\in\{0,1\}^8\subset\F_p^8.
\]
Thus \(b_i\) records the neighbors of the \(i\)-th remaining vertex among
the first eight vertices. Since \(\rk M=\rk B\), the Schur
complement of \(B\) vanishes:
\begin{equation}\label{eq:d6-schur}
        E=N^TRN.
\end{equation}
The image of \(M\) is
\[
        \operatorname{im}_{\F_p}M=
        \left\{
        \begin{pmatrix}v\\N^TRv\end{pmatrix}:v\in\F_p^8
        \right\}.
\]
By Lemma~\ref{lem:d6-kernel}, the system \(Mw=\one_{28}\) is solvable over
\(K\). Since its coefficients lie in \(\F_p\), it is already solvable over
\(\F_p\). Hence
\begin{equation}\label{eq:d6-ones-image}
        N^TR\one_8=\one_{20}.
\end{equation}
Consequently,
\[
        z=\begin{pmatrix}R\one_8\\\mathbf0_{20}\end{pmatrix}
        \in\F_p^{28}
\]
satisfies \(Mz=\one_{28}\). Let \(y\) be the zero-sum solution from
Lemma~\ref{lem:d6-kernel}. After extending scalars to \(K\), we have
\(y-z\in\ker M=\ker\widehat H\), so \(y-z\) also has coordinate sum zero.
Therefore
\begin{equation}\label{eq:d6-ones-isotropic}
        \one_8^TR\one_8
        =\one_{28}^Tz
        =\one_{28}^Ty
        =0.
\end{equation}

The diagonal and off-diagonal entries of equation~\eqref{eq:d6-schur}
give equations~\eqref{eq:d6-norm} and~\eqref{eq:d6-compat}, respectively.
Equations~\eqref{eq:d6-ones-image} and~\eqref{eq:d6-ones-isotropic} give
equations~\eqref{eq:d6-main} and~\eqref{eq:d6-isotropic}. Finally, the \(b_i\) are
distinct: if \(b_i=b_j\) for \(i\ne j\), then
\[
        E_{ij}=b_i^TRb_j=b_i^TRb_i=-\lambda,
\]
contrary to \(E_{ij}\in\{0,1\}\) and \(\lambda\ne0,-1\).
\end{proof}

\subsection{Exact verification of the finite problem}\label{subsec:d6-computation}
The preceding reduction involves all odd primes. The next proposition first
reduces the problem to finitely many characteristics and then checks the
remaining cases by exact computation.

A self-contained SageMath program implementing this verification is available
in the
\href{https://github.com/kyleyip111/two-distance-sets-over-finite-fields}
{GitHub repository for this paper}.
It uses SageMath's interface to \texttt{nauty} to generate all graphs on eight
vertices up to isomorphism, performs the resultant and finite-field
computations, and uses the exact \texttt{Cliquer} routine to compute the clique
numbers of the final auxiliary graphs. No floating-point arithmetic is used.

\begin{proposition}\label{prop:d6-computer}
There do not exist an odd prime \(p\), a scalar
\(\lambda\in\F_p\setminus\{0,-1\}\), a simple graph \(\Gamma\) on eight
vertices, and $20$ binary column vectors satisfying
equations~\eqref{eq:d6-norm}--\eqref{eq:d6-isotropic},
where \(C\in M_8(\F_p)\) is the reduction modulo \(p\) of the
\(0\)-\(1\) adjacency matrix of \(\Gamma\) and
\(C-\lambda I_8\) is nonsingular.
\end{proposition}

\begin{proof}
We first reduce the problem to a finite computation.  For an eight-vertex graph
\(\Gamma\), let \(C_\Gamma\in M_8(\mathbb Z)\) be its \(0\)-\(1\)
adjacency matrix, and define polynomials over \(\mathbb Z\) by
\[
        \Delta_\Gamma(t)=\det(C_\Gamma-tI_8),
        \qquad
        S_\Gamma(t)=\one_8^T\adj(C_\Gamma-tI_8)\one_8.
\]
For \(b\in\{0,1\}^8\subset\mathbb Z^8\), put
\[
        F_{\Gamma,b}(t)=b^T\adj(C_\Gamma-tI_8)b+t\Delta_\Gamma(t),
\]
\[
        G_{\Gamma,b}(t)=b^T\adj(C_\Gamma-tI_8)\one_8-\Delta_\Gamma(t).
\]
Fix an odd prime \(p\), let \(\overline C_\Gamma\in M_8(\F_p)\) be the
reduction of \(C_\Gamma\), and let \(\lambda\in\F_p\) satisfy
\(\Delta_\Gamma(\lambda)\ne0\).  Since
\[
        (\overline C_\Gamma-\lambda I_8)^{-1}
        =\Delta_\Gamma(\lambda)^{-1}
          \adj(\overline C_\Gamma-\lambda I_8),
\]
equations~\eqref{eq:d6-norm}, \eqref{eq:d6-main}, and
\eqref{eq:d6-isotropic}, with \(C=\overline C_\Gamma\), are equivalent to
\begin{equation}\label{eq:d6-polynomials}
        F_{\Gamma,b}(\lambda)=G_{\Gamma,b}(\lambda)
        =S_\Gamma(\lambda)=0
        \qquad\text{in }\F_p.
\end{equation}

Write
\[
        S_\Gamma(t)=c\prod_{j=1}^r h_j(t)^{e_j}
\]
in \(\mathbb Z[t]\), where the \(h_j\) are primitive and irreducible.
The content \(c\) divides the leading coefficient \(-8\), so \(p\nmid c\).
Consequently, if \(S_\Gamma(\lambda)=0\) in \(\F_p\), then
\(h_j(\lambda)=0\) for at least one \(j\).  Fix such a factor \(h=h_j\).  Since
\(\lambda\ne0,-1\) and \(\Delta_\Gamma(\lambda)\ne0\), we may ignore the
factors \(t\), \(t+1\), and those dividing \(\Delta_\Gamma\).

For each remaining factor \(h\), define
\[
        \mathcal B_h=
        \{b\in\{0,1\}^8:h\mid F_{\Gamma,b}\text{ and }
        h\mid G_{\Gamma,b}\}.
\]
For \(b\notin\mathcal B_h\), put
\[
        r_{h,b}=\gcd\bigl(\Res_t(h,F_{\Gamma,b}),
        \Res_t(h,G_{\Gamma,b})\bigr)\in\mathbb Z_{>0}.
\]
Here the subscript \(t\) indicates that the resultants eliminate the only
indeterminate \(t\), so both resultants are integers.  Since \(h\) is irreducible over
\(\mathbb Q\), the equality \(\Res_t(h,F_{\Gamma,b})=0\) is equivalent to
\(h\mid F_{\Gamma,b}\), and similarly for \(G_{\Gamma,b}\).  Thus the two
resultants cannot both vanish when \(b\notin\mathcal B_h\), and hence
\(r_{h,b}>0\).

Suppose that the system of equations~\eqref{eq:d6-polynomials} holds modulo
an odd prime \(p\) at a root \(\lambda\) of \(h\).  Then the reductions of
\(h\) and \(F_{\Gamma,b}\), and also those of \(h\) and
\(G_{\Gamma,b}\), have the common root \(\lambda\).  Hence both resultants
vanish modulo \(p\), so either \(b\in\mathcal B_h\) or \(p\mid r_{h,b}\).
Thus there can be $20$ vectors satisfying equations~\eqref{eq:d6-norm}
and~\eqref{eq:d6-main} only if
\begin{equation}\label{eq:d6-exceptional-filter}
        |\mathcal B_h|+
        \#\{b\notin\mathcal B_h:p\mid r_{h,b}\}\ge20.
\end{equation}

The remaining calculation is exact. SageMath's \texttt{nauty} interface
generates one representative of each of the \(12{,}346\) isomorphism classes
of graphs on eight vertices. For every such graph, the program processes each
relevant factor \(h\) and all \(256\) binary vectors. It finds
\(13{,}809\) relevant pairs \((\Gamma,h)\) and verifies
\[
        |\mathcal B_h|\le16
\]
in every case. Consequently, if
inequality~\eqref{eq:d6-exceptional-filter} holds, then \(p\mid r_{h,b}\) for
at least one \(b\notin\mathcal B_h\). Since only finitely many positive
integers \(r_{h,b}\) occur, this reduces the problem to finitely many primes.
The program factors these integers exactly; no upper bound for \(p\) is
imposed in advance. Only $60$ odd primes remain, the largest of which is
\(1367\).

For each retained graph, factor, and prime, the program finds the roots
\(\lambda\in\F_p\) of \(h\), removes \(0\) and \(-1\), and then removes
duplicate triples. This gives \(17{,}347\) triples
\((\Gamma,p,\lambda)\). Of these, \(6{,}480\) satisfy
\(\Delta_\Gamma(\lambda)=0\) and are discarded, leaving \(10{,}867\)
nonsingular triples. For each remaining triple, the program tests all
\(256\) binary vectors in \(\F_p^8\). Exactly \(4{,}075\) triples admit at
least $20$ vectors satisfying equations~\eqref{eq:d6-norm} and
\eqref{eq:d6-main}, and they occur only in the following characteristics:
\begin{center}
\begin{tabular}{@{}rr@{}}
\toprule
Characteristic & Number of triples\\
\midrule
\(3\)  & \(2455\)\\
\(5\)  & \(1562\)\\
\(7\)  & \(54\)\\
\(19\) & \(2\)\\
\(47\) & \(2\)\\
\bottomrule
\end{tabular}
\end{center}
For each such triple, form an auxiliary graph
\(\mathcal C_{\Gamma,p,\lambda}\). Its vertices are the binary vectors
\(b\in\{0,1\}^8\) satisfying equations~\eqref{eq:d6-norm} and
\eqref{eq:d6-main}. Equation~\eqref{eq:d6-compat} gives the pairwise
compatibility condition: two distinct vertices \(b,c\) are adjacent precisely
when
\[
        b^T(\overline C_\Gamma-\lambda I_8)^{-1}c\in\{0,1\}
        \qquad\text{in }\F_p.
\]
Here equation~\eqref{eq:d6-isotropic}, which depends only on
\((\Gamma,p,\lambda)\), has already been imposed. The graph \(\Gamma\)
describes the eight vertices corresponding to the chosen principal
submatrix, whereas the vertices of
\(\mathcal C_{\Gamma,p,\lambda}\) represent the possible neighborhood
vectors of the remaining vertices of the original graph. Consequently, the
$20$ vectors produced by a hypothetical \(28\)-point configuration would
form a clique of order \(20\) in
\(\mathcal C_{\Gamma,p,\lambda}\). For each of the \(4{,}075\) auxiliary
graphs, the program computes the clique number exactly
using SageMath's \texttt{Cliquer} interface. The largest clique number obtained
is \(17\), so none contains a clique of order \(20\).
\end{proof}

\begin{proof}[Proof of Theorem~\ref{thm:blokhuis}\textup{(ii)}]
Proposition~\ref{prop:d6-construction} gives a \(27\)-point example. If a
set of size at least \(28\) existed, choose a \(28\)-point subset in
which both distances still occur. Proposition~\ref{prop:d6-reduction} would
then produce a system excluded by Proposition~\ref{prop:d6-computer}.
Therefore every such set has at most \(27\) points.
\end{proof}

\section{Concluding remarks}
The finite-field situation separates two facts that coincide over the real numbers. The
polynomial method gives the general Larman--Rogers--Seidel bound
\(|S|\le \binom{d+2}{2}+d+1\), while Blokhuis' sharper real bound requires a
positivity step. Over finite fields, this step is replaced by the rank defects
in Theorem~\ref{thm:refined-rank-sensitive}.

Theorem~\ref{thm:blokhuis} determines exactly when the Blokhuis bound is
attained. The midpoint construction and the Johnson--coset construction,
together with the characteristic \(3\) examples, attain it in every dimension \(d\ne6\).
Dimension $6$ is exceptional: no standard space over an odd finite field
contains \(28\) such points, and the exact maximum is \(27\). On the other
hand, the characteristic \(3\) construction shows that the Blokhuis bound is
not a universal finite-field upper bound; the full Larman--Rogers--Seidel bound
is attained in infinitely many dimensions.

\section*{Acknowledgments}
The authors thank Petr Lison\v{e}k, Shuxing Li, Hong-Jun Ge, and Wei-Hsuan Yu for helpful discussions. The research of the first author was supported in part by an
NSERC Discovery Grant and by the National Research, Development, and Innovation Office of
Hungary, NKFIH, Grant No. KKP133819 and Excellence 151341.

\end{document}